\documentclass[11pt]{amsart}
\usepackage[all]{xy}
\usepackage{amsmath}
\usepackage{amsfonts}
\usepackage{amssymb}
\usepackage{latexsym}
\usepackage{graphicx}
\usepackage{color}
\usepackage{amscd}
\usepackage{epsfig}
\usepackage{cite}

\newtheorem{problem}{Problem}[section]
\newtheorem{definition}[problem]{Definition}
\newtheorem{lemma}[problem]{Lemma}
\newtheorem{theorem}[problem]{Theorem}

\newtheorem{corollary}[problem]{Corollary}

\newtheorem{conjecture}[problem]{Conjecture}

\title{$T$-adic exponential sums of polynomials in one variable}
\author{Chunlei Liu}
\address{Department of Mathematics, Shanghai Jiao Tong
University, Shanghai 200240, P.R. China, E-mail: clliu@sjtu.edu.cn}
\author{Wenxin Liu}\address{School of Mathematical Sciences, Beijing Normal
University, Beijing 100875, P.R. China, E-mail:
wenxin8210@mail.bnu.edu.cn}
\thanks{This research is supported by NSFC Grant No.
10671015.}
\begin{document}
\maketitle
\begin{abstract}
The $T$-adic exponential sum of a polynomial in one variable is
studied. An explicit arithmetic polygon in terms of the highest two
exponents of the polynomial is proved to be a lower bound of the
Newton polygon of the $C$-function of the T-adic exponential sum.
This bound gives lower bounds for the Newton polygon of the
$L$-function of exponential sums of $p$-power order.
\end{abstract}



\section{Introduction}

Let $p$ be a prime number,  $q$ a power of $p$, and $\mathbb{F}_q$
the finite field with $q$ elements. Let $W$ be the Witt ring scheme,
$\mathbb{Z}_{q}=W(\mathbb{F}_q)$, and
$\mathbb{Q}_{q}=\mathbb{Z}_{q}[\frac{1}{p}]$ .

Let $\triangle\supsetneq\{0\}$ be an integral convex polytope in
$\mathbb{R}^n$, and $I$ the set of vertices of $\triangle$ different
from the origin. Let
$$f(x)=\sum\limits_{u\in\triangle}(a_ux^u,0,0,\cdots)\in W(\mathbb{F}_q[x_1^{\pm1},x_2^{\pm1},\cdots,x_n^{\pm1}])\text{ with }
\prod_{u\in I}a_u\neq0,$$ where $x^u=x_1^{u_1}x_2^{u_2}\cdots
x_n^{u_n}$ if $u=(u_1,u_2,\cdots,u_n)\in \mathbb{Z}^n$.

\begin{definition}For any positive integer $l,$
the $T$-adic exponential sum associated to $f$ is the sum:
$$S_{f}(l,T)=\sum\limits_{x\in (\mathbb{F}_{q^l}^{\times})^n}
(1+T)^{Tr_{\mathbb{Q}_{q^l}/\mathbb{Q}_p}(f(x))}\in\mathbb{Z}_p[[T]].$$
The $T$-adic $L$-function $L_f(s,T)$ associated to $f$ is defined by
the formula
$$L_f(s,T)=\exp(\sum\limits_{l=1}^{\infty}S_f(l,T)\frac{s^l}{l})\in 1+s\mathbb{Z}_p[[T]][[s]].$$
\end{definition}
Let $m\geq1$, $\zeta_{p^m}$ a primitive $p^m$-th root of unity, and
$\pi_m=\zeta_{p^m}-1$. Then the specialization
 $L_{f}(s,\pi_m)$ is the $L$-function of
$p$-power order exponential sums $S_{f}(l,\pi_m)$. These $p$-power
order exponential sums were studied by Adolphson-Sperber \cite{AS}
for $m=1$, and by Liu-Wei \cite{LW} for $m\geq1$.

We view $L_f(s, T)$ as a power series in the single variable $s$
with coefficients in the $T$-adic  complete field
$\mathbb{Q}_p((T))$. So we call it a $T$-adic $L$-function.

\begin{definition}
The $T$-adic $C$-function $C_f(s,T)$ associated to $f$ is the
generating function
$$C_f(s,T) =\exp(\sum\limits_{l=1}^{\infty}-(q^l-1)^{-n}S_{f}(l,T)\frac{s^l}{l}).$$
\end{definition}

We have
$$L_{f}(s,T) = \prod_{i=0}^n C_{f}(q^is,T)^{(-1)^{n-i+1}{n\choose i}},$$
and $$C_{f}(s,T)^{(-1)^{n-1}}= \prod_{j=0}^{\infty}
L_{f}(q^js,T)^{n+j-1\choose j}.$$ So the $C$-function $C_f(s,T)$ and
the $L$-function $L_f(s,T)$ determine each other. From the last
identity, one sees that
$$C_f(s,T)\in
1+s\mathbb{Z}_p[[T]][[s]].$$We also view $C_f(s, T)$ as a power
series in the single variable $s$ with coefficients in the $T$-adic
complete field $\mathbb{Q}_p((T))$. The $C$-function $C_f(s,T)$ was
shown to be $T$-adic entire by Liu-Wan \cite{LWn}.
\\

Let $C(\triangle)$ be the  cone generated by $\triangle$, and
$M(\triangle)=C(\triangle)\cap \mathbb{Z}^n$. There is a degree
function $\deg$ on $C(\triangle)$ which is
$\mathbb{R}_{\geq0}$-linear and takes the values $1$ on every
co-dimension $1$ face not containing $0$.  For example, if
$\triangle\subset\mathbb{Z}$, then the degree function on
$C(\triangle)$ is defined by the formula
$$\deg(a)=\left\{
                                   \begin{array}{ll}
                                     0, & \hbox{} a=0,\\
                                     \frac{a}{d({\rm sgn}(a))}, &
\hbox{}a\neq0,
                                   \end{array}
                                 \right.
$$
where $d(\varepsilon)$ is the nonzero endpoint of $\triangle$ with
sign $\varepsilon=\pm1$. For $a\not\in C(\triangle)$, we define
$\deg(a)=+\infty$.
\begin{definition}A convex function on $\mathbb{R}_{\geq0}$ which is
linear between consecutive integers with initial value $0$ is called
the Hodge polygon of $\triangle$ if its slopes between consecutive
integers are the numbers $\deg(a)$, $a\in M(\triangle)$. We denote
it by  $H_{\triangle}^\infty$.\end{definition}

Liu-Wan \cite{LWn}
 proved the following.
\begin{theorem}We have $$T-adic \text{ NP of } C_{f}(s,T)\geq {\rm ord}_p(q)(p-1)H_{\triangle}^\infty,$$
 where NP is the short for Newton polygon.
 \end{theorem}

From now on we assume that $\triangle=[0,d]$.
\begin{definition}For $a\in \mathbb{N}$, we define
$$\delta_{\in}(a)=\left\{
                           \begin{array}{ll}
                            1, & \hbox{ }pi\equiv a(d)\text{ for some }i<d\{\frac{a}{d}\}, \\
                            0, & \hbox{ otherwise,}
                            \end{array}
                            \right.
$$where $\{\cdot\}$ is the fractional part of a real number.
\end{definition}

\begin{definition}A convex function on
$\mathbb{R}_{\geq 0}$ which is linear between consecutive integers
with initial value $0$ is called the arithmetic polygon of
$\triangle =[0,d]$ if  its slopes between consecutive integers are
the numbers
$$\varpi_{\triangle}(a)=\lceil(p-1)\deg(a)\rceil-\delta_{\in}(a), a\in \mathbb{N},$$
where $\lceil\cdot\rceil$ is the least integer equal or greater than
a real number.  We denote it by $p_{\triangle}$.
\end{definition}

Liu-Liu-Niu \cite{LLN} proved the following.
\begin{theorem}If $p>3d$,  then
$$T-adic \text{ NP of } C_{f}(s,T)\geq \text{ord}_p(q)p_{\triangle},
$$with equality holding for a generic $f$ of degree $d$.
\end{theorem}

By a result of Li\cite{Li}, $L_{f}(s,\pi_m)$ is a polynomial with
degree $p^{m-1}d$ if $p\nmid d$. Combined this result with the above
theorem, one can infer the following.
\begin{theorem}If $p>3d$,  then
$$\pi_m-adic \text{ NP of } L_{f}(s,\pi_m)\geq \text{ord}_p(q)p_{\triangle}\text{ on }[0,p^{m-1}d],$$
with equality holding for a generic $f$ of degree $d$.
\end{theorem}
The Newton polygon of the $L$-function $L_{f}(s,\pi_m)$  for $m=1$
was studied by Zhu \cite{Zhu, Zhu2} and Blache-F\'{e}rard \cite{BF}.

We assume that the second highest exponent of $f$ is $ k$. So $k
\leq d-1$, and
$$f(x)=(a_dx^d,0,0,\cdots )+\sum\limits_{i=1}^k(a_ix^i,0,0,\cdots )\in W(\mathbb{F}_q[x])\text{ with }a_da_k\neq
0.$$

For  $a\in  \mathbb{N}$, define

\begin{align*}
\varpi_{d,[0,k]}(a)=&[\frac{pa}{d}]-[\frac{a}{d}]+[\frac{r_{pa}}{k}]-[\frac{r_a}{k}]+\sum\limits_{i=1}^{r_a}(
1_{\{\frac{r_{pi}}{k}\}>\{\frac{r_a}{k}\}}-1_{\{\frac{r_i}{k}\}>\{\frac{r_a}{k}\}})\\
&-\sum\limits_{i=1}^{r_{a-1}}(
1_{\{\frac{r_{pi}}{k}\}>\{\frac{r_{a-1}}{k}\}}-
1_{\{\frac{r_i}{k}\}>\{\frac{r_{a-1}}{k}\}}),
\end{align*}
where  $r_a=d\{\frac{a}{d}\}$ for $a\in \mathbb{N}$.

\begin{definition}A convex function on
$\mathbb{R}_{\geq 0}$ which is linear between consecutive integers
with initial value $0$ is called the arithmetic polygon of
$\{d\}\cup[0,k]$ if its slopes between consecutive integers are the
numbers $\varpi_{d,[0,k]}(a)$, $a\in \mathbb{N}$. We denote it by
$p_{d,[0,k]}$ .
\end{definition}

We can prove the following.
\begin{theorem} \label{upper-bound} We have $p_{d,[0,k]}\geq p_{\triangle}$.
\end{theorem}
 The  main result of this paper is
the following.
\begin{theorem}\label{first}If $p>d(2d+1)$,  then $$T-adic \text{ NP of } C_{f}(s,T)\geq
ord_p(q)p_{d,[0,k]}.$$
\end{theorem}
\begin{corollary}\label{second} If $p>d(2d+1)$, then
  $$\pi_m\text{ -adic NP of }
L_{f}(s,\pi_m) \geq \text{ord}_p(q)p_{d,[0,k]}  \text{ on }
[0,p^{m-1}d].$$
\end{corollary}

We put forward the following conjecture.
\begin{conjecture}Let $k<d$ be fixed positive integers. If $p$ is sufficiently large, then for a generic $f$ of the form
$$f(x)=(a_dx^d,0,0,\cdots )+\sum\limits_{i=1}^k(a_ix^i,0,0,\cdots )\in W(\mathbb{F}_q[x])\text{ with }a_da_k\neq
0,$$ we have $$T-adic \text{ NP of } C_{f}(s,T)=
ord_p(q)p_{d,[0,k]}.$$\end{conjecture} The above conjecture implies
the following.
\begin{conjecture}Let $k<d$ be fixed positive
integers. Then for a generic $f$ of the form
$$f(x)=(a_dx^d,0,0,\cdots )+\sum\limits_{i=1}^k(a_ix^i,0,0,\cdots )\in W(\mathbb{F}_q[x])\text{ with }a_da_k\neq
0,$$ the $T^{p-1}$-adic Newton polygon of $C_{f}(s,T)$ uniformly
converges to ${\rm ord}_p(q)H_{[0,d]}^{\infty}$ as $p$ goes to
infinity.\end{conjecture} The above conjecture again implies the
following.
\begin{conjecture}Let $k<d$ be fixed positive
integers. Then for a generic $f$ of the form
$$f(x)=(a_dx^d,0,0,\cdots )+\sum\limits_{i=1}^k(a_ix^i,0,0,\cdots )\in W(\mathbb{F}_q[x])\text{ with }a_da_k\neq
0,$$ the $\pi_m^{p-1}$-adic Newton polygon of $C_{f}(s,\pi_m)$
converges to ${\rm ord}_p(q)H_{[0,d]}^{\infty}$ as $p$ goes to
infinity.\end{conjecture} The last two conjectures are analogues of
a conjecture put forward by Wan \cite{Wa,Wa2}.

\section{The $T$-adic Dwork Theory}In this section we review the $T$-adic analogue of Dwork
theory on exponential sums.

Let
$$E(t)=\exp(\sum_{i=0}^{\infty}\frac{t^{p^i}}{p^i})=\sum\limits_{i=0}^{+\infty}\lambda_it^i
\in 1+t{\mathbb Z}_p[[t]]$$ be the $p$-adic Artin-Hasse exponential
series. Define a new $T$-adic uniformizer $\pi$ of ${\mathbb
Q}_p((T))$ by the formula $E(\pi)=1+T$. Let $\pi^{1/d}$ be a fixed
$d$-th root of $\pi$. Let $a\mapsto\hat{a}$ be the Teichm\"{u}ller
lifting. One can show that the series
$$E_f(x) :=E(\pi \hat{a}_dx^d)\prod\limits_{i=1}^{k}E(\pi \hat{a}_ix^i)$$ lies in the $T$-adic
Banach module
$$L=\{\sum_{i\in
\mathbb{N}}c_i\pi^{\deg(i)}x^i:\
 c_i\in\mathbb{Z}_q[[\pi^{1/d}]] \}.$$ Note
that the Galois group of $\mathbb{Q}_q$ over $\mathbb{Q}_p$ can act
on $L$ but keeping $\pi^{1/d}$ as well as the variable $x$ fixed.
Let $\sigma$ be the Frobenius element in the Galois group such that
$\sigma(\zeta)=\zeta^p$ if $\zeta$ is a $(q-1)$-th root of unity.
  Let $\Psi_p$ be the operator on $L$ defined by the formula
$$\Psi_p(
 \sum\limits_{i\in
\mathbb{N}} c_ix^i)=\sum\limits_{i\in \mathbb{N}} c_{pi}x^i.$$ Then
$\Psi:=\sigma^{-1}\circ\Psi_p\circ E_f$ acts on the $T$-adic Banach
module
$$B=\{\sum\limits_{i\in \mathbb{N}}c_i\pi^{\deg(i)}x^i \in L,\
 \text{\rm ord}_T(c_i)\rightarrow+\infty
 \text{ if }\deg(i)\rightarrow+\infty\}.$$
We call it Dwork's $T$-adic semi-linear operator because it is
semi-linear over $\mathbb{Z}_q[[\pi^{\frac{1}{d}}]].$  Let $b=
\text{\rm ord}_p(q)$, the $b$-iterate $\Psi^b$ is linear over
$\mathbb{Z}_q[[\pi^{1/d}]]$, since
$$\Psi^{b}=\Psi_p^{b}\circ
\prod\limits_{i=0}^{b-1}E_{f}^{\sigma^i}(x^{p^i}).$$ One can show
that $\Psi$ is completely continuous in the sense of Serre
\cite{Se}. So $\det(1-\Psi^bs\mid
B/\mathbb{Z}_{q}[[\pi^{\frac{1}{d}}]])$ and $\det(1-\Psi s\mid
B/\mathbb{Z}_p[[\pi^{\frac{1}{d}}]])$ are well-defined.

We now state the $T$-adic Dwork trace formula.
\begin{theorem}[$T$-adic Dwork trace formula \cite{LWn}]\label{analytic-trace-formula}
We have
$$C_f(s,T)=\det(1-\Psi^bs\mid
B/\mathbb{Z}_{q}[[\pi^{\frac{1}{d}}]]).$$
\end{theorem}
\section{Key estimate}
In order to study $$C_f(s,T)=\det(1-\Psi^bs\mid
B/\mathbb{Z}_{q}[[\pi^{\frac{1}{d}}]]),$$ we first study
$$\det(1-\Psi s\mid B/\mathbb{Z}_p[[\pi^{\frac{1}{d}}]])=\sum\limits_{i=0}^{+\infty}(-1)^ic_is^i.$$
We are going to show that
\begin{theorem}  \label{finalreduction}If $p>d(2d+1),$ then we
have $${\rm ord}_{\pi}(c_{bm})\geq bp_{d,[0,k]}(m).$$
\end{theorem}

Fix a normal basis $\bar{\xi}_1,\cdots,\bar{\xi}_b$ of
$\mathbb{F}_q$ over $\mathbb{F}_p$. Let $\xi_1,\cdots,\xi_b$ be
their Teichm\"{u}ller lifts. Then $\xi_1,\cdots,\xi_b$ is a normal
basis of $\mathbb{Q}_q$ over $\mathbb{Q}_p$, and $\sigma$ acts on
$\xi_1,\cdots,\xi_b$ as a permutation. Let
$(\gamma_{(i,u),(j,\omega)})_{i,j\in \mathbb{N},1\leq u,\omega\leq
b}$ be the matrix of $\Psi$ on
$B\otimes_{\mathbb{Z}_p}\mathbb{Q}_p(\pi^{1/d})$ with respect to the
basis $\{\xi_ux^i\}_{i\in \mathbb{N},1\leq u\leq b}$. Then
$$c_{bm}=\sum\limits_{R}\text{det}((\gamma_{(i,u),(j,\omega)})_{(i,u),(j,\omega)\in
R}),
$$
where $R$ runs over all subsets of $
\mathbb{N}\times\{1,2,\cdots,b\}$ with cardinality $bm$.

So Theorem \ref{finalreduction} is reduced to the following.
\begin{theorem}\label{det}
Let $R\subset \mathbb{N}\times\{1,2,\cdots,b\}$ be a subset of
cardinality $bm$. If $p>d(2d+1)$, then
$${\rm ord}_T(\det(\gamma_{(i,u),(j,\omega)})_{(i,u),(j,\omega)\in
R})\geq bp_{d,[0,k]}(m).$$
\end{theorem}

Let $O(\pi^{\alpha})$ denote any element of $\pi$-adic order $\geq
\alpha$.
\begin{lemma}We have $$\gamma_{(i,u),(j,\omega)}=O(\pi^{[\frac{pi-j}{d}]+\lceil\frac{r_{pi-j}}{k}\rceil}).$$
\end{lemma}
\proof Write $$E_f(x)= \sum\limits_{i\in \mathbb{N}}\gamma_ix^i.$$
 Then $$\gamma_i=\sum\limits_{\stackrel{dn_d+\sum\limits_{j=1}^kjn_j=i}
{n_j\geq0}}\pi^{\sum\limits_{j=1}^kn_j+n_d}
\prod\limits_{j=1}^k\lambda_{n_j}\hat{a}_j^{n_j}\lambda_{n_d}\hat{a}_d^{n_d}=O(\pi^{[\frac{i}{d}]+\lceil\frac{r_i}{k}\rceil}).$$
Since
$$(\xi_{\omega}\gamma_{pi-j})^{\sigma^{-1}}=\sum\limits_{u=1}^{b
}\gamma_{(i,u),(j,\omega)}\xi_{u},$$ we have $${\rm
ord}_T(\gamma_{(i,u),(j,\omega)})={\rm ord}_T(\gamma_{pi-j}).$$The
lemma now follows.\endproof

By the above lemma, Theorem \ref{det} is reduced to the following.
\begin{theorem}
Let $R\subset \mathbb{N}\times\{1,2,\cdots,b\}$ be a subset of
cardinality $bm$, and $\tau$ a permutation of $R$. Suppose that
$p>d(2d+1)$. Then
$$\sum\limits_{(i,u)\in
R}([\frac{pi-\tau(i)}{d}]+\lceil\frac{r_{pi-\tau(i)}}{k}\rceil)\geq
bp_{d,[0,k]}(m),$$ where $\tau(i)$ is defined by
$\tau(i,u)=(\tau(i), \tau(u))$.
\end{theorem}
\proof  By definition, we have
\begin{align*}
 p_{d,[0,k]}(m)&=\sum\limits_{a=1}^{m-1}\varpi_{d,[0,k]}(a)\\
&=\sum\limits_{a=1}^{m-1}([\frac{pa}{d}]-[\frac{a}{d}]+[\frac{r_{pa}}{k}]-[\frac{r_a}{k}]+
1_{\{\frac{r_{pa}}{k}\}>\{\frac{r_{m-1 }}{k}\}}-
1_{\{\frac{r_a}{k}\}>\{\frac{r_{m -1}}{k}\}}).
\end{align*}
Note that \begin{align*} &\sum\limits_{(i,u)\in
R}([\frac{pi-\tau(i)}{d}]+\lceil\frac{r_{pi-\tau(i)}}{k}\rceil)\\
&=\sum\limits_{(i,u)\in
R}([\frac{pi}{d}]-[\frac{i}{d}]+[\{\frac{pi}{d}\}-\{\frac{\tau(i)}{d}\}]+
\lceil\frac{r_{pi}-r_{\tau(i)}-d[\{\frac{pi}{d}\}-\{\frac{\tau(i)}{d}\}]}{k}\rceil)\\
&=\sum\limits_{(i,u)\in
R}([\frac{pi}{d}]-[\frac{i}{d}]+\lceil\frac{r_{pi}-r_{\tau(i)}+(d-k)1_{r_{\tau(i)}>r_{pi}}}{k}\rceil)\\
&=\sum\limits_{(i,u)\in
R}([\frac{pi}{d}]-[\frac{i}{d}]+[\frac{r_{pi}}{k}]-[\frac{r_i}{k}]+\lceil\frac{(d-k)1_{r_{\tau(i)}>r_{pi}}}{k}
+\{\frac{r_{pi}}{k}\}-\{\frac{r_{\tau(i)}}{k}\}\rceil).
\end{align*}

We  have
\begin{align*}
& \sum\limits_{(i,u)\in
R}\lceil\frac{(d-k)1_{r_{\tau(i)}>r_{pi}}}{k}
+\{\frac{r_{pi}}{k}\}-\{\frac{r_{\tau(i)}}{k}\}\rceil\\
&\geq \sum\limits_{(i,u)\in R} 1_{\{\frac{r_{\tau(i) }}{k}\}\leq
\{\frac{r_{m-1}
}{k}\}<\{\frac{r_{pi}}{k}\}}\\
&\geq \sum\limits_{(i,u)\in R }
1_{\{\frac{r_{pi}}{k}\}>\{\frac{r_{m-1}}{k}\}}-
\sum\limits_{(i,u)\in R} 1_{\{\frac{r_i}{k}\}>\{\frac{r_{m-1}
}{k}\}}.
\end{align*}
We also have
\begin{align*}
&\sum\limits_{(i,u)\in
R}([\frac{pi}{d}]-[\frac{i}{d}]+[\frac{r_{pi}}{k}]-[\frac{r_i}{k}])\\
&=b\sum\limits_{i=1}^{m-1}([\frac{pi}{d}]-[\frac{i}{d}]+[\frac{r_{pi}}{k}]-[\frac{r_i}{k}])
+\sum\limits_{\stackrel{(i,u)\in R}{i\geq
m}}([\frac{pi}{d}]-[\frac{i}{d}]+[\frac{r_{pi}}{k}]-[\frac{r_i}{k}])\\
&\quad -\sum\limits_{\stackrel{(i,u)\notin R}{0\leq
i<m}}([\frac{pi}{d}]-[\frac{i}{d}]+[\frac{r_{pi}}{k}]-[\frac{r_i}{k}])\\
\end{align*}
\begin{align*}
&\geq b\sum\limits_{i=1}^{m-1}([\frac{pi}{d}]-[\frac{i}{d}]+[\frac{r_{pi}}{k}]-[\frac{r_i}{k}])+N([\frac{pm}{d}]-[\frac{m}{d}]-[\frac{d-1}{k}])\\
&-N([\frac{p(m-1)}{d}]-[\frac{m-1}{d}]+[\frac{d-1}{k}])\\
&\geq
b\sum\limits_{i=1}^{m-1}([\frac{pi}{d}]-[\frac{i}{d}]+[\frac{r_{pi}}{k}]-[\frac{r_i}{k}])+N([\frac{p}{d}]-1-2(d-1)).
\end{align*}
where $N=\#\{(i,u)\in R| i\geq m\}=\#\{(i,u)\notin R|0\leq i <m\}$.

Similarly, we have $$\sum\limits_{(i,u)\in
R}1_{\{\frac{r_{pi}}{k}\}>\{\frac{r_{m-1}}{k}\}}\geq
b\sum\limits_{i=1}^{m-1}1_{\{\frac{r_{pi}}{k}\}>\{\frac{r_{m-1}}{k}\}}-N,
$$
and $$\sum\limits_{(i,u)\in R}
1_{\{\frac{r_i}{k}\}>\{\frac{r_{m-1}}{k}\}}\leq
 b\sum\limits_{i=1}^{m-1}1_{\{\frac{r_i}{k}\}>\{\frac{r_{m-1}}{k}\}}+N.
$$
 Therefore
\begin{align*}
&\sum\limits_{(i,u)\in R}([\frac{pi-\tau(i)}{d}]+\lceil\frac{r_{pi-\tau(i)}}{k}\rceil)\\
&\geq
b\sum\limits_{i=1}^{m-1}([\frac{pi}{d}]-[\frac{i}{d}]+[\frac{r_{pi}}{k}]-[\frac{r_i}{k}])+b\sum\limits_{i=1}^{m-1}
1_{\{\frac{r_{pi}}{k}\}>\{\frac{r_{m-1}}{k}\}}
-b\sum\limits_{i=1}^{m-1}1_{\{\frac{r_i}{k}\}>\{\frac{r_{m-1}}{k}\}}\\
&\quad+N\big([\frac{p}{d}]-1-2(d-1)-2\big)\\
&= bp_{d,[0,k]}(m)+N\big([\frac{p}{d}]-2d+1\big) \geq
bp_{d,[0,k]}(m).
\end{align*}\endproof

\section{Proof of the main result}
In this section we prove Theorem \ref{first} , which says that, if
$p>d(2d+1)$,  then $$T-adic \text{ NP of } C_{f}(s,T)\geq
ord_p(q)p_{d,[0,k]}.$$

\begin{lemma}\label{q2p}The Newton polygon of
 $\det(1-\Psi^bs^b\mid
B/\mathbb{Z}_{q}[[\pi^{\frac{1}{d}}]])$ coincides with that of
$\det(1-\Psi s\mid B/\mathbb{Z}_p[[\pi^{\frac{1}{d}}]])$.\end{lemma}
\proof Note that
$$\det(1-\Psi^b s\mid B/\mathbb{Z}_p[[\pi^{\frac{1}{d}}]])
={\rm Norm}(\det(1-\Psi^b s\mid
B/\mathbb{Z}_q[[\pi^{\frac{1}{d}}]])),$$ where Norm is the norm map
from $\mathbb{Z}_q[[\pi^{\frac{1}{d}}]]$ to
$\mathbb{Z}_p[[\pi^{\frac{1}{d}}]]$. The lemma now follows from the
equality
$$\prod\limits_{\zeta^b=1}\det(1-\Psi\zeta s\mid B/\mathbb{Z}_p[[\pi^{\frac{1}{d}}]])
=\det(1-\Psi^b s^b\mid B/\mathbb{Z}_p[[\pi^{\frac{1}{d}}]]).$$ \qed

\begin{corollary}\label{coefficient-q2p}
The $T$-adic Newton polygon of
 $\det(1-\Psi^bs\mid
B/\mathbb{Z}_{q}[[\pi^{\frac{1}{d}}]])$ is the lower convex closure
of the points
$$(m,\text{ord}_{T}(c_{bm})),\ m=0,1,\cdots.$$\end{corollary}

\proof By Lemma \ref{q2p}, the $T$-adic Newton polygon of
 $\det(1-\Psi^bs^b\mid
B/\mathbb{Z}_{q}[[\pi^{\frac{1}{d}}]])$ is the lower convex closure
of the points
$$(i,\text{ord}_{T}(c_i)),\ i=0,1,\cdots.$$
It is clear that $(i,\text{ord}_{T}(c_i))$ is not a vertex of that
polygon if $b\nmid i$. So that Newton polygon is the lower convex
closure of the points
$$(bm,\text{ord}_{T}(c_{bm})),\ m=0,1,\cdots.$$
It follows that the $T$-adic Newton polygon of
 $\det(1-\Psi^bs\mid
B/\mathbb{Z}_{q}[[\pi^{\frac{1}{d}}]])$ is the lower convex closure
of the points
$$(m,\text{ord}_{T}(c_{bm})),\ m=0,1,\cdots.$$
\endproof

We now prove  Theorem \ref{first}.

{\it Proof of Theorem \ref{first}. } By Theorem
\ref{analytic-trace-formula}, we have
$$C_f(s,T)=\det(1-\Psi^bs\mid
B/\mathbb{Z}_{q}[[\pi^{\frac{1}{d}}]]).$$ Then by Corollary
\ref{coefficient-q2p}, the $T$-adic Newton polygon of $C_f(s,T)$ is
the lower convex closure of the points
$$(m,\text{ord}_{T}(c_{bm})),\ m=0,1,\cdots.$$
Therefore the result follows from Theorem \ref{finalreduction},
which says that, if $p>d(2d+1),$ then we have
$${\rm ord}_{\pi}(c_{bm})\geq bp_{d,[0,k]}(m).$$\qed

We conclude this section by proving Corollary \ref{second}.

 {\it
Proof of Corollary \ref{second}. } Assume that $ L_f(s,
\pi_m)=\prod\limits_{i=1}^{p^{m-1}d}(1-\beta_is) $. Then
$$C_f(s,\pi_m )= \prod\limits_{j=0}^{\infty} L_f(q^js, \pi_m)=
\prod\limits_{j=0}^{\infty}\prod\limits_{i=1}^{p^{m-1}d}(1-\beta_iq^js).$$
Therefore the slopes of the  $q$-adic Newton polygon of
$C_{f}(s,\pi_m)$ are the numbers
$$j+{\rm ord}_q(\beta_i),\
1\leq i\leq p^{m-1}d, j=0,1,\cdots.$$ It is well-known that ${\rm
ord}_q(\beta_i)\leq 1$ for all $i$. Therefore, $$q-adic \text{ NP of
}L_{f}(s,\pi_m)= q-adic\text{ NP of }C_{f}(s,\pi_m) \text { on }
[0,p^{m-1}d].$$ It follows that
 $$\pi_m-adic \text{ NP of
}L_{f}(s,\pi_m)= \pi_m-adic\text{ NP of }C_{f}(s,\pi_m) \text { on }
[0,p^{m-1}d].$$ By  the integrality of $C_f(s,T)$ and Theorem
\ref{first}, we have
$$\pi_m-adic \text{ NP of } C_{f}(s,\pi_m)\geq T-adic \text{ NP of } C_{f}(s,T)\geq \text{ord}_p(q)p_{d,[0,k]}.$$

Therefore,   $$\pi_m\text{ -adic NP of } L_{f}(s,\pi_m) \geq
\text{ord}_p(q)p_{d,[0,k]}  \text{ on } [0,p^{m-1}d].$$ \qed

\section{comparison between arithmetic polygons}
In this section we prove Theorem \ref{upper-bound}, which says that
$$p_{d,[0,k]}\geq p_{\triangle}.$$

{\it Proof of Theorem \ref{upper-bound} }  It is clear that
$p_{d,[0,k]}(0)=p_{\triangle}(0)$.  It suffices to show that, for
$m\in \mathbb{N}$, we have $p_{d,[0,k]}(m+1) \geq
p_{\triangle}(m+1)$. By a result in Liu-Liu-Niu \cite{LLN}, we have

 $$
p_{\triangle}(m+1)=\sum\limits_{a=1}^{m }\varpi_{\triangle}(a)=
\sum\limits_{a=1}^{m }(\lceil \frac{pa}{d}\rceil-\lceil
\frac{a}{d}\rceil)+\sum\limits_{a=1}^{r_{m }}1_{r_{pa}>r_{m }}.
 $$
 By definition, we have
\begin{align*}
 p_{d,[0,k]}(m+1)&=\sum\limits_{a=1}^{m}\varpi_{d,[0,k]}(a)\\
&=\sum\limits_{a=1}^{m }([\frac{pa}{d}]-[\frac{a}{d}])+
\sum\limits_{a=1}^{r_{m }}([\frac{r_{pa}}{k}]-[\frac{a}{k}]+
1_{\{\frac{r_{pa}}{k}\}>\{\frac{r_{m }}{k}\}}-
1_{\{\frac{a}{k}\}>\{\frac{r_{m }}{k}\}}),
\end{align*}

For $m\geq 0$,
$$\sum\limits_{a=1}^{m }([\frac{pa}{d}]-[\frac{a}{d}])=\sum\limits_{a=1}^{m }(\lceil
\frac{pa}{d}\rceil-\lceil \frac{a}{d}\rceil),$$

Let  $A_1=\{1\leq a\leq r_m|a\neq r_{pi} \text{ for some } 1\leq
i\leq r_m\}$. Note that $$\{a:1 \leq a\leq r_m\} =A_1\cup A_2,$$
 where $A_2= \{r_{pa}|1 \leq a, r_{pa}\leq r_m\}.$
And
$$\{r_{pa}|1\leq a\leq r_m\}=A_2\cup A_3,$$
where $A_3=\{r_{pa}>r_m| 1 \leq a\leq r_m\}$, so we have
$|A_1|=|A_3|$. Then

\begin{align*}
&\sum\limits_{a=1}^{r_m}([\frac{r_{pa}}{k}]-[\frac{a}{k}]+
 1_{\{\frac{r_{pa}}{k}\}>\{\frac{r_m}{k}\}}- 1_{\{\frac{a}{k}\}>\{\frac{r_m}{k}\}})\\
&=\sum\limits_{r_{pa}\in A_3}[\frac{r_{pa}}{k}]-\sum\limits_{a\in
A_1}[\frac{a}{k}]+ \sum\limits_{r_{pa}\in
A_3}1_{\{\frac{r_{pa}}{k}\}>\{\frac{r_m}{k}\}}- \sum\limits_{a\in
A_1} 1_{\{\frac{a}{k}\}>\{\frac{r_m
}{k}\}}\\
&=\sum\limits_{r_{pa}\in
A_3}([\frac{r_{pa}}{k}]-[\frac{r_m}{k}]+1_{\{\frac{r_{pa}}{k}\}>\{\frac{r_m}{k}\}})+
\sum\limits_{a\in A_1}([\frac{r_m
}{k}]-[\frac{a}{k}]-1_{\{\frac{a}{k}\}>\{\frac{r_m
}{k}\}})\\
&=\sum\limits_{r_{pa}\in A_3 }\lceil\frac{r_{pa}-r_m}{k}\rceil+
\sum\limits_{a\in
A_1}(\lceil\frac{r_m-a}{k}\rceil-1_{\{\frac{a}{k}\}>\{\frac{r_m}{k}\}}-1_{\{\frac{a}{k}\}<\{\frac{r_m
}{k}\}})\\
&\geq |A_3|=\sum\limits_{a=1}^{r_m}1_{r_{pa}>r_m}.
\end{align*}

The  theorem now follows. \qed

\end{document}